\documentclass[12pt]{article}

\usepackage{amsmath,amsfonts}
\usepackage{graphicx}
\usepackage{multirow}
\newtheorem{thm}{Theorem}
\newtheorem{alg}{Algorithm}
\newtheorem{ass}{Assumption}
\newcommand{\R}{\mathbb{R}}

\title{Convex inner approximations of nonconvex semialgebraic sets
applied to fixed-order controller design\footnote{A preliminary version
of this work was presented during the International Symposium on Mathematical
Theory of Networks and Systems, Budapest, Hungary, 5-9 July 2010.}}

\begin{document}

\author{Didier Henrion$^{1,2,3}$, Christophe Louembet$^{1,2}$}

\footnotetext[1]{CNRS; LAAS; 7 avenue du colonel Roche, F-31077 Toulouse; France.}
\footnotetext[2]{Universit\'e de Toulouse; UPS, INSA, INP, ISAE; UT1, UTM, LAAS; F-31077 Toulouse; France}
\footnotetext[3]{Faculty of Electrical Engineering, Czech Technical University in Prague,
Technick\'a 2, CZ-16626 Prague, Czech Republic}

\maketitle

\begin{abstract}
We describe an elementary algorithm to build convex inner approximations of nonconvex sets.
Both input and output sets are basic semialgebraic sets given as lists of defining multivariate polynomials.
Even though no optimality guarantees can be given (e.g. in terms of volume maximization for bounded sets),
the algorithm is designed to preserve convex boundaries as much as possible, while removing regions
with concave boundaries. In particular, the algorithm leaves invariant a given convex set.
The algorithm is based on Gloptipoly 3, a public-domain Matlab package solving nonconvex polynomial optimization
problems with the help of convex semidefinite programming (optimization over linear matrix inequalities,
or LMIs). We illustrate how the algorithm can be used to design fixed-order controllers for linear systems,
following a polynomial approach.
\end{abstract}


\section{Introduction}

The set of controllers stabilizing a linear system is generally {\em nonconvex}
in the parameter space, and this is an essential difficulty faced by
numerical algorithms of computer-aided control system design, see
e.g. \cite{hsk03} and references therein. It follows from the derivation
of the Routh-Hurwitz stability criterion (or its discrete-time counterpart)
that the set of stabilizing controllers is real basic semialgebraic, i.e.
it is the intersection of sublevel sets of given multivariate polynomials.
A convex inner approximation of this nonconvex semialgebraic stability
region was obtained in \cite{hsk03} in the form of linear matrix inequalities
(LMI) obtained from univariate polynomial positivity conditions, see also
\cite{kkl07}. Convex polytopic
inner approximations were also obtained in \cite{n06}, for discrete-time
stability, using reflection coefficients. Convex inner approximations make it
possible to design stabilizing controllers with the help
of convex optimization techniques, at the price of loosing
optimality w.r.t. closed-loop performance criteria ($H_2$ norm,
$H_{\infty}$ norm or alike).

Generally speaking, the technical literature abounds of convex {\em outer}
approximations of nonconvex semialgebraic sets. In particular, such
approximations form the basis of many branch-and-bound global optimization
algorithms \cite{n04}. By construction, Lasserre's hierarchy of LMI relaxations
for polynomial programming is a sequence of embedded convex outer
approximations which are semidefinite representable, i.e. which are
obtained by projecting affine sections of the
convex cone of positive semidefinite matrices, at the price of
introducing lifting variables \cite{hl04}.

After some literature search, we could not locate any systematic
constructive procedure to generate convex {\em inner} approximations
of nonconvex semialgebraic sets, contrasting sharply with the many
convex outer approximations mentioned above. In the context of fixed-order
controller design, inner approximations correspond to a guarantee of
stability, at the price of loosing optimality. No such stability
guarantee can be ensured with outer approximations.

The main contribution of this paper is therefore an elementary algorithm,
readily implementable in Matlab, that generates convex inner approximations
of nonconvex sets.
Both input and output sets are basic semialgebraic sets given as lists
of defining multivariate polynomials. Even though no optimality guarantees
can be given in terms of volume maximization for bounded sets,
the algorithm is designed to preserve convex boundaries as much as possible,
while removing regions with concave boundaries. In particular,
the algorithm leaves invariant a given convex set.
The algorithm is based on Gloptipoly 3, a public-domain Matlab package
solving nonconvex polynomial optimization problems with the help of
convex LMIs \cite{g3}. Even though the algorithm can be useful
on its own, e.g. for testing convexity of semialgebraic sets,
we illustrate how it can be used to design
fixed-order controllers for linear systems,
following a polynomial approach.

\section{Convex inner approximation}

Given a basic closed semialgebraic set
\begin{equation}\label{levelset_S}
S=\{x\in\R^n:p_1(x)\leq 0\:\ldots\:p_m(x)\leq 0\}
\end{equation}
where $p_i$ are multivariate polynomials,
we are interested in computing another basic
closed semialgebraic set
\begin{equation}\label{levelset_Sinner}
\bar{S}=\{x\in\R^n:\bar{p}_1(x)\leq 0\:\ldots\:\bar{p}_{\bar{m}}(x)\leq 0\}
\end{equation}
which is a valid convex inner approximation of $S$, in the sense that
\[
\bar{S} \subset S.
\]
Ideally, we would like to find the tightest possible approximation,
in the sense that the complement set $S \backslash \bar{S} = \{x \in S \: :\:
x \notin \bar{S}\}$ is as small as
possible. Mathematically we may formulate the problem as the
volume minimization problem
\[
\inf_{\bar{S}} \int_{S \backslash \bar{S}} dx
\]
but since set $S$ is not necessarily bounded we should make sure
that this integral makes sense. Moreover, computing
the volume of a given semialgebraic set is a difficult task in general
\cite{sirev}, so we expect that optimizing such a quantity is as much
as difficult. In practice, in this paper, we will
content ourselves of an inner approximation that removes the nonconvex
parts of the boundary and keeps the convex parts as much as possible.

\section{Detecting nonconvexity}

Before describing the method,
let us recall some basics definitions on polynomials and
differential geometry.
Let $x \in \R^n \mapsto p_i(x)\in \R[x]$ be a multivariate polynomial
of total degree $d$.
Let
\[
g_i(x)=\left[\frac{\partial p_i(x)}{\partial x_j}\right]_{j=1\ldots n} \in \R^n[x]
\]
be its gradient vector and
\[
H_i(x) = \left[\frac{\partial^2 p_i(x)}{\partial x_j \partial x_k}\right]_{j,k=1\ldots n} \in \R^{n\times n}[x]
\]
its (symmetric) Hessian polynomial matrix.
Define the optimization problem
\begin{equation}\label{opti}
\begin{array}{llll}
q_i & = & \min_{x,y} & y^TH_i(x,y)y \\
&& \mathrm{s.t.} & p_i(x) = 0 \\
&& & p_j(x) \leq 0,\: j=1\ldots m,\: j \neq i \\
&& & y^T g_i(x) = 0 \\
&& & y^T y = 1
\end{array}
\end{equation}
with global minimizers $\{x^1\ldots x^{k_i}\}$
and $\{y^1\ldots y^{k_i}\}$.

Let us make the following nondegeneracy assumption on
defining polynomials $p_i(x)$:
\begin{ass}\label{sing}
There is no point $x$ such that $p_i(x)$ and $g_i(x)$
vanish simultaneously while satisfying $p_j(x) \leq 0$
for $j=1,\ldots,m$, $j \neq i$.
\end{ass}
Since the polynomial system $p_i(x) = 0$, $g_i(x) = 0$,
involves $n+1$ equations for $n$ unknowns,
Assumption \ref{sing} is satisfied generically. In other words,
in the Euclidean space of coefficients of polynomials $p_i(x)$,
instances violating Assumption \ref{sing} belong to a variety
of Lebesgue measure zero, and an arbitrarily small perturbation
on the coefficients generates a perturbed set
$S_{\epsilon}$ satisfying Assumption \ref{sing}.

\begin{thm}\label{convexity_thm}
Under Assumption \ref{sing},
polynomial level set \eqref{levelset_S}
is convex if and only if $q_i \geq 0$ for all $i=1,\ldots,m$.
\end{thm}

{\bf Proof:}
The boundary of set $S$ consists of points $x$ such that
$p_i(x) = 0$ for some $i$, and $p_j(x) \leq 0$ for $j \neq i$.
In the neighborhood of such a point, consider the Taylor series
\begin{equation}\label{taylor}
p_i(x+y) = p_i(x) + y^T g_i(x) + y^T H_i(x) y + O(y^3)
\end{equation}
where $O(y^3)$ denotes terms of degree $3$ or higher
in entries of vector $y$, the local coordinates.
By Assumption \ref{sing}, the gradient $g_i(x)$ does not
vanish along the boundary, and
hence convexity of the boundary is inferred from the
quadratic term in expression (\ref{taylor}). More specifically,
when $y^T g_i(x) = 0$, vector $y$ belongs to the
hyperplane tangent to $S$ at point $x$. Let $V$ be a matrix
spanning this linear subspace of dimension $n-1$
so that $y = V \hat{y}$ for some $\hat{y}$. The quadratic form
$y^T H_i(x) y = \hat{y}^T V^T H_i V \hat{y}$ can be diagonalised
with the congruence transformation $\hat{y} = U \bar{y}$
(Schur decomposition), and hence
$y^T H_i(x) y = \bar{y}^T U^T V^T H_i V U \bar{y}^T = \sum_{i=1}^{n-1}
h_i(x) \bar{y}^2_i$. The eigenvalues $h_i(x)$, $i=1,\ldots,n-1$
are reciprocals of the principal curvatures of the surface. Problem (\ref{opti})
then amounts to finding the minimum curvature, which is non-negative when
the surface is locally convex around $x$.$\Box$

In the case of three-dimensional surfaces ($n=3$), the ideas of
tangent plane, local coordinates and principal curvatures used
in the proof of Theorem \ref{convexity_thm} are
standard notions of differential geometry, see e.g. Section 3.3. in \cite{docarmo}
for connections between principal curvatures
and eigenvalues of the local Hessian form (called the second fundamental form,
once suitably normalized).

\begin{figure}[htb]
\centering{\includegraphics[width=0.6\columnwidth]{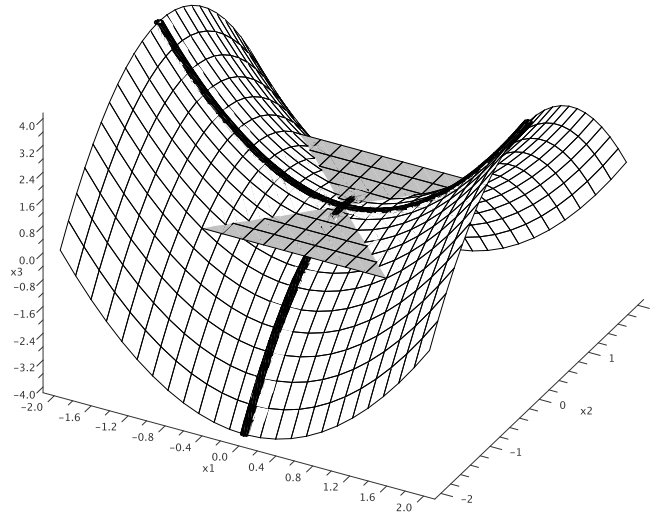}}
\caption{Hyperboloid of one sheet (white), with tangent plane (gray)
at the origin, a saddle point with a tangent convex parabola (thick black)
and a tangent concave hyperbola (thick black).\label{hyperboloid}}
\end{figure}

As an example illustrating the proof of Theorem \ref{convexity_thm},
consider the hyperboloid of one sheet
$S = \{x \in {\mathbb R}^3 \: :\: p_1(x) = x_1^2 - x_2^2 - x_3 \leq 0\}$
with gradient and Hessian
\[
g_1(x) = \left[\begin{array}{c}2x_1\\-2x_2\\-1\end{array}\right], \quad
H_1(x) = \left[\begin{array}{ccc}2&0&0\\0&-2&0\\0&0&0\end{array}\right].
\]
At the origin $x=0$, the tangent plane is $T = \{y \in {\mathbb R}^3
\: :\: y_3 = 0\}$ and $p_1(y) = 2y_1^2-2y_2^2$ is a bivariate
quadratic form with eigenvalues $2$ and $-2$, corresponding respectively
to the convex parabola $\{x \: :\: x_2^2+x_3 = 0\}$
(positive curvature) and concave hyperbola $\{x \: :\: x_1^2-x_3 = 0\}$
(negative curvature),
see Figure \ref{hyperboloid}.

Theorem \ref{convexity_thm} can be exploited in an algorithmic way
to generate a convex inner approximation of a semialgebraic set.

\begin{alg}\label{inner_alg}{\bf (Convex inner approximation)}
\begin{description}
\item {\sl Input:} Polynomials $p_i$, $i=1\ldots m$ defining set $S$
as in (\ref{levelset_S}). Small nonnegative scalar $\epsilon$.
\item {\sl Output:} Polynomials $\bar{p}_i$, $i=1\ldots\bar{m}$
defining set $\bar{S}$ as in (\ref{levelset_Sinner}).
\item {\sl Step 1:} Let $i=1$.
\item {\sl Step 2:} If $\deg p_i \leq 1$ then go to Step 5.
\item {\sl Step 3:}  If $p_i(x)\in S$, solve optimization problem (\ref{opti})
for optimum $q_i$ and minimizers $\{x^1\ldots x^k\}$. If $p_i(x)\notin S$, go to Step 5.
\item {\sl Step 4:} If $q_i < 0$, then select one of the minimizers $x^j$,
$j=1\ldots k_i$, let $p_{m+1} = g_i(x^j)(x-x^j)+\epsilon$.
Then let $m=m+1$, and go to step 3.
\item {\sl Step 5:} Let $i=i+1$. If $i \leq m$ then go to Step 2.
\item {\sl Step 6:} Return $\bar{p}_i=p_i$, $i=1,\ldots m$.
\end{description}
\end{alg}

The idea behind the algorithm is as follows.
At Step 3, by solving the polynomial optimization problem
of Theorem \ref{convexity_thm}
we identify a point of minimal curvature along
algebraic varieties defining the boundary of $S$.
If the minimal curvature is negative, then
we separate the point from the set with
a gradient hyperplane, and we iterate on
the resulting semialgebraic set.
At the end, we obtain a valid inner approximation.

Note that Step 2 checks if the boundary is affine,
in which case the minimum curvature is zero and there is
no optimization problem to be solved.

The key parameter of the algorithm is the small
positive scalar $\epsilon$ used at Step 4 for separating
strictly a point of minimal curvature, so that
the algorithm does not identify it again
at the next iteration.
Moreover, in Step 4, one must elect arbitrarily a minimizer. We will discuss this issue later in this paper.

Finally, as pointed out to us by a referee, the ordering of
the sequence of input polynomials $p_i$ has an impact on the sequence
of output polynomials $\bar{p}_i$, and especially on the size
of the convex inner approximation $\bar{S}$. However, it seems
very difficult to design a priori an optimal ordering policy.

\section{Matlab code and geometric examples}

At each step of Algorithm \ref{inner_alg} we have to solve a potentially
nonconvex polynomial optimization problem. For that purpose, we
use Gloptipoly 3, a public-domain Matlab package \cite{g3}. The
methodology consists in building and solving a hierarchy of
embedded linear matrix inequality (LMI) relaxations of the polynomial
optimization problem, see the survey \cite{l09}.
The LMI problems are solved numerically with
the help of any semidefinite programming solver (by default
Gloptipoly 3 uses SeDuMi). Under the assumption that our original
semi-algebraic set is compact, the sequence of minimizers obtained
by solving the LMI relaxations is ensured to converge mononotically
to the global minimum. Under the additional assumption that the global
optima live on a zero-dimensional variety (i.e. there is a finite
number of them), Gloptipoly 3 eventually extracts some of them
(not necessarily all, but at least one) using numerical linear
algebra. The LMI problems in the
hierarchy have a growing number of variables and constraints,
and the main issue is that we cannot predict in advance how
large has to be the LMI problem to guarantee global optimality.
In practice however we observe that it is not necessary to
go very deep in the hierarchy to have a numerical certificate
of global optimality.

\subsection{Hyperbola}

Let us first with the elementary example of an unbounded nonconvex
hyperbolic region $S = \{x \in {\mathbb R}^2 \: :\: p_1(x) \leq 0\}$
with $p_1(x)=-1+x_1x_2$, for which optimization problem (\ref{opti})
reads
\[
\begin{array}{ll}
\min & 2y_1y_2 \\
\mathrm{s.t.} & x_2 y_1 + x_1 y_2 = 0 \\
& -1 + x_1 x_2 = 0 \\
& y_1^2+y_2^2 = 1. \\
\end{array}
\]
Necessary optimality conditions yield immediately
$k_1=2$ global minimizers $x^1 = \frac{\sqrt{2}}{2}(1,1)$,
$y^1 = \frac{\sqrt{2}}{2}(1,-1)$
and $x^2 = \frac{\sqrt{2}}{2}(-1,-1)$, $y^2 = \frac{\sqrt{2}}{2}(-1,1)$,
and hence two additional (normalized) affine constraints
$p_2(x) = -2 + x_1 + x_2$ and $p_3(x) = -2 - x_1 - x_2$
defining the slab $\bar{S} = \{x \: :\: p_i(x) \leq 0,
\: i=1,2,3\} = \{x \: :\: -2 \leq x_1 + x_2 \leq 2\}$
which is indeed a valid inner approximation of $S$.

\subsection{Egg quartic}

\begin{figure}[htb]
\centering{\includegraphics[width=0.75\columnwidth]{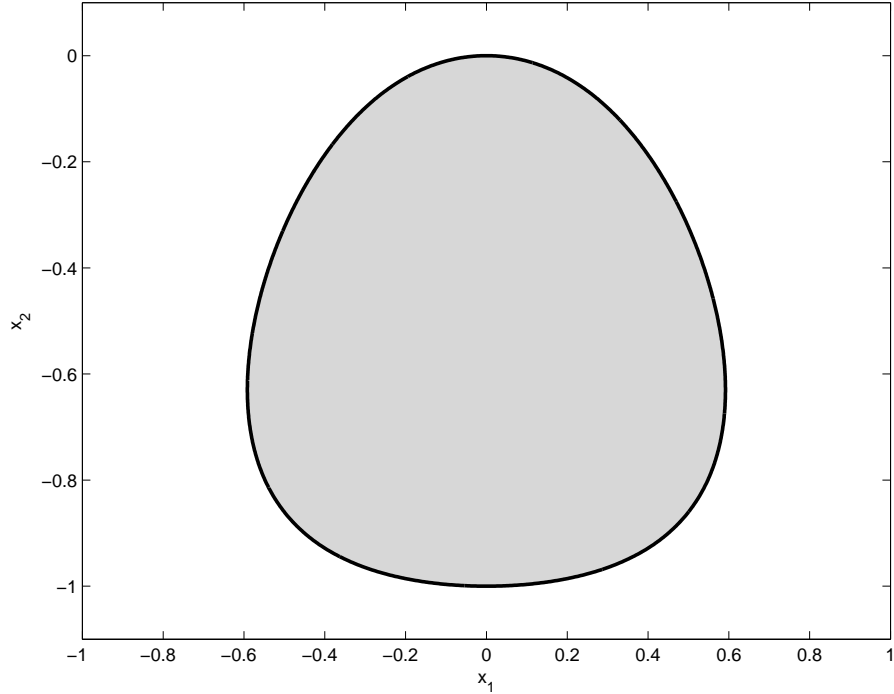}}
\caption{Convex smooth quartic.
\label{convexquartic}}
\end{figure}

Now we show that Algorithm \ref{inner_alg} can be used to detect convexity
of a semialgebraic set. Consider the smooth quartic sublevel set
$S = \{x \in {\mathbb R}^2 \: :\: p_1(x)=x_1^4+x_2^4+x_1^2+x_2 \leq 0\}$
represented on Figure \ref{convexquartic}.
Assumption \ref{sing} is ensured since
the gradient $g_1(x)=[2x_1(x_1^2+2) \:\: 4x_2^3+1]$ cannot
vanish for real $x$.

A Matlab implementation of the first steps of the
algorithm can be easily written using Gloptipoly 3:
\begin{verbatim}
% problem data
mpol x y 2
p1 = x(1)^4+x(2)^4+x(1)^2+x(2);
g1 = diff(p1,x); % gradient
H1 = diff(g1,x); % Hessian
% LMI relaxation order
order = 3;
% build LMI relaxation
P = msdp(min(y'*H1*y), p1==0, ...
    g1*y==0, y'*y==1, order);
% solve LMI relaxation
[status,obj] = msol(P)
\end{verbatim}
Notice that we solve the LMI relaxation of order 3 (e.g. moments
of degree 6) of problem (\ref{opti}). In GloptiPoly, an LMI relaxation
is solved with the command {\tt msol} which returns two output arguments:
{\tt status} and {\tt obj}. Argument {\tt status} can take the following values:
\begin{itemize}
\item {\tt -1} if the LMI relaxation is infeasible or could not be solved
for numerical reasons;
\item {\tt 0} if the LMI relaxation could be solved but it is impossible
to detect global optimality and to extract global optimizers, in which case
{\tt obj} is a lower (resp. upper) bound on the global minimum (resp. maximum)
of the original optimization problem;
\item {\tt +1} if the LMI relaxation could be solved, global optimality
is certified and global minimizers are extracted, in which case {\tt obj}
is the global optimum of the original optimization problem.
\end{itemize}
Running the above script, Gloptipoly returns {\tt obj = 2.0000}
and {\tt status = 1}, certifying that the minimal curvature
is strictly positive, and hence that the polynomial sublevel set
is convex.

Note that in this simple case, convexity of set $S$ follows
directly from positive semidefiniteness of the Hessian
$H_1(x) = \mathrm{diag}\:(12x_1^2+2,\:12x_2^2)$,
yet Algorithm \ref{inner_alg} can systematically
detect convexity in more complicated cases.

\subsection{Waterdrop quartic}\label{secwaterdrop}

Consider the quartic $S = \{x \in {\mathbb R}^2 \: :\:
p_1(x) = x_1^4+x_2^4+x_1^2+x_2^3 \leq 0\}$ which has a singular point
at the origin, hence violating Assumption \ref{sing}.

\begin{figure}[htb]
\centering{\includegraphics[width=0.75\columnwidth]{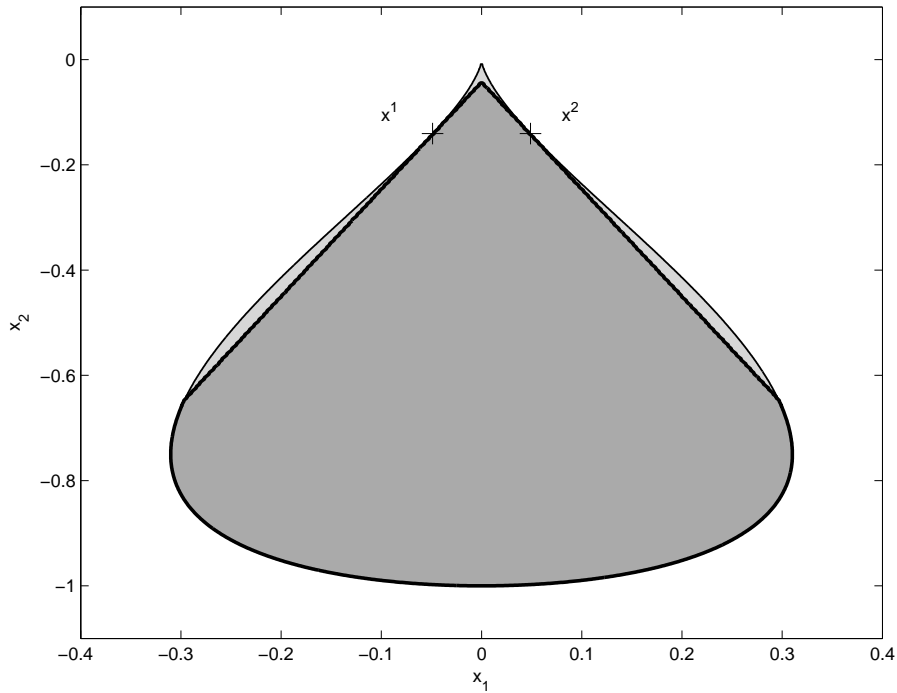}}
\caption{Nonconvex waterdrop quartic (light gray) and its
convex inner approximation (dark gray) obtained by adding
affine constraints at two points $x^1$ and $x^2$
of minimal curvature.\label{waterdrop}}
\end{figure}

Applying Algorithm \ref{inner_alg},
the LMI relaxation of order 4 (moments of degree 8)
yields a globally minimal curvature of $-0.094159$
achieved at the 2 points
$x^1=(-0.048892,\:-0.14076)$ and $x^2=(0.048896,\:-0.14076)$.
With the two additional affine constraints
$p_k(x) = g_1(x^k)(x-x^k) \leq 0$, $k=2,3$, the resulting set
$\bar{S}$ has a globally
minimal curvature of $1$ certified at the LMI relaxation of order 4,
and therefore it is a valid convex inner approximation of $S$, see
Figure \ref{waterdrop}.

This example illustrates that Algorithm \ref{inner_alg} can
work even when Assumption \ref{sing} is violated. Here the
singularity is removed by the additional affine constraints.
This example also shows that symmetry of the problem
can be exploited, since two global minimizers are
found (distinct points with the same minimal curvature)
to remove two nonconvex parts of the boundary simultaneously.

\subsection{Singular quartic}\label{ss_singular_quartic}

Consider the quartic $S = \{x \in {\mathbb R}^2 \: :\:
p_1(x) = x_1^4+x_2^4+x_2^3 \leq 0\}$ which has a singular point
at the origin, hence violating Assumption \ref{sing}.

\begin{figure}[htb]
\centering{\includegraphics[width=0.75\columnwidth]{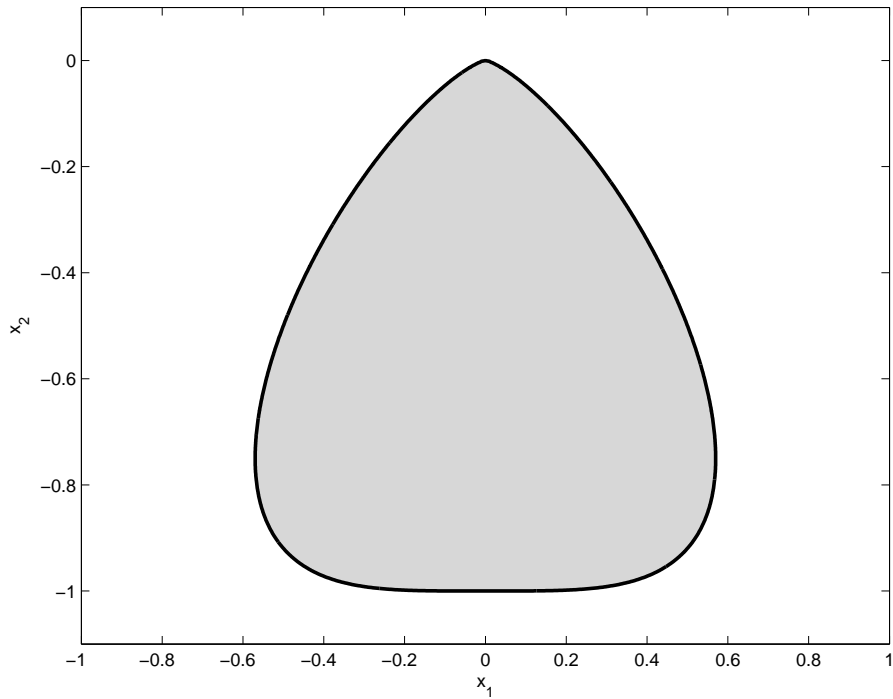}}
\caption{Numerically convex singular quartic.\label{singquartic}}
\end{figure}

Running Algorithm \ref{inner_alg}, we obtain the following sequence
of bounds on the minimum curvature, for increasing LMI relaxation orders:
\[
\begin{array}{l|cccc}
{\tt order} & 2 & 3 & 4 & 5 \\\hline
{\tt obj} & -7.5000\cdot10^{-1} & -7.7502\cdot10^{-2} &
-8.5855\cdot10^{-3} & -4.9525\cdot10^{-3}
\end{array}
\]
GloptiPoly is not able to certify global optimality, so we can
only speculate that the global minimum is zero and hence that
set $S$ is convex, see Figure \ref{singquartic}. We may say
that set $S$ is numerically convex.

\begin{figure}[htb]
\centering{\includegraphics[width=0.75\columnwidth]{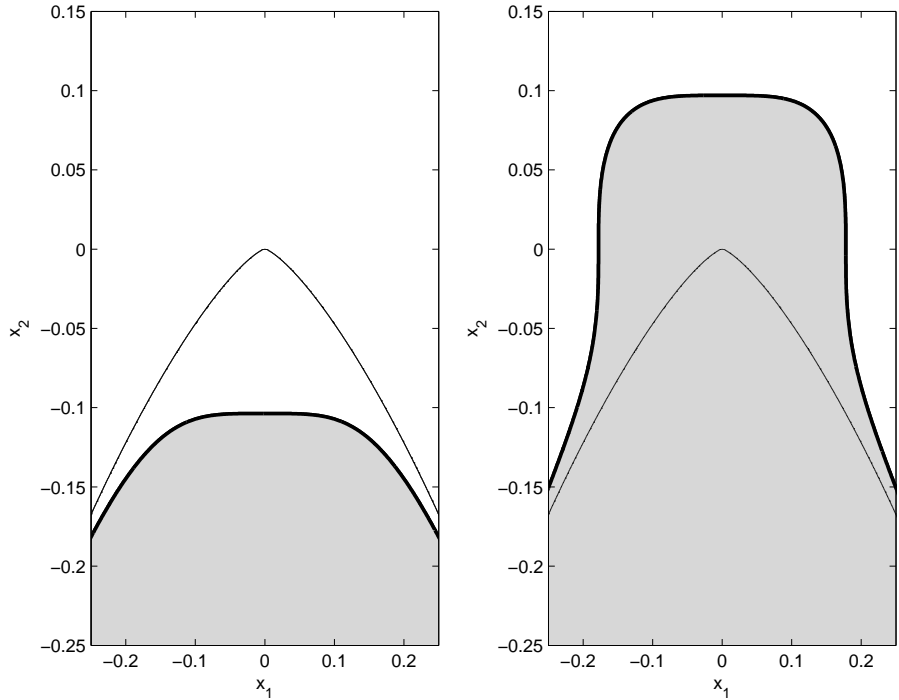}}
\caption{Perturbed quartic $p_1(x)+\epsilon \leq 0$ (bold line) can be convex ($\epsilon=10^{-3}$) or nonconvex ($\epsilon=-10^{-3}$) near
singularity of original quartic level set $p_1(x) = 0$ (light line).
\label{singquarticzoom}}
\end{figure}

Indeed if we strenghten the constraint $p_1(x) \leq 0$ into
$p_1(x)+\epsilon \leq 0$ for a small positive $\epsilon$, say $10^{-3}$, then
GloptiPoly 3 certifies global optimality and convexity
with {\tt obj = -4.0627e-7}
at the 4th LMI relaxation. On the other hand, if we relax the
constraint into $p_1(x)+\epsilon \leq 0$ with
a negative $\epsilon=-10^{-3}$, then GloptiPoly 3 certifies global
optimality and nonconvexity with {\tt obj = -0.22313}
at the 4th LMI relaxation. We can conclude that the optimum
of problem \ref{opti} is sensitive, or ill-conditioned,
with respect to the problem data, the coefficients of $p_1(x)$.
The reason behind this ill-conditioning is the singularity
of $S$ at the origin, see Figure \ref{singquarticzoom}
which represents the effect of perturbing the constraint
$p_1(x) \leq 0$ around the singularity.

\section{Control applications}

In this section we focus on control applications of Algorithm \ref{inner_alg}, which is used to generate convex inner approximation
of stability regions in the parameter space.

\subsection{Third-order discrete-time stability region}

Algorithm \ref{inner_alg} can lend insight into the (nonconvex) geometry
of the stability region. Consider the simplest non-trivial case of
a third-order discrete-time polynomial $x_1+x_2z+x_3z^2+z^3$
which is stable (roots within the open unit disk) if and only if
parameter $x=(x_1,x_2,x_3)$ lies within the interior of
compact region $S = \{x \in {\mathbb R}^3 \: :\:
p_1(x) = -x_1-x_2-x_3-1\leq 0,\: p_2(x) = x_1-x_2+x_3-1\leq 0,\:
p_3(x) = x_1^2-x_1x_3+x_2-1\leq 0\}$. Stability region $S$ is
nonconvex, delimited by two planes $p_1(x) = 0$, $p_2(x) = 0$
and a hyperbolic paraboloid $p_3(x) = 0$
see e.g. \cite[Example 11.4]{ackermann}.

Optimization problem (\ref{opti}) corresponding to convexity check
of the hyperbolic paraboloid reads as follows:
\begin{equation}\label{schur3_opti}
\begin{array}{ll}
\min & -2y_1^2+2y_1y_3 \\
\mathrm{s.t.} & x_1^2-x_1x_3+x_2-1 = 0 \\
& -x_1-x_2-x_3-1 \leq 0 \\
& x_1-x_2+x_3-1 \leq 0 \\
& (2x_1-x_3)y_1 +y_2 +x_3y_3 = 0 \\
& y_1^2+y_2^2+y_3^2 = 1.
\end{array}
\end{equation}
The objective function and the last constraint depend only on $y$,
and necessary optimality conditions obtained by differentiating
the Lagrangian $-2y_1^2+2y_1y_3+t(y_1^2+y_2^2+y_3^2-1)$ with
respect to $y$ yield the symmetric pencil equation
\[
\left[\begin{array}{lll}
-4+2t & 0 & 2 \\
0 & 2t & 0 \\
2 & 0 & 2t
\end{array}\right]
\left[\begin{array}{c}
y_1 \\ y_2 \\ y_3
\end{array}\right] = 0
\]
From the determinant
of the above 3-by-3 matrix, equal to $t(t^2-2t-1)$, we
conclude that multiplier $t$ can be equal to $1-\sqrt{2}$,$0$
or $1+\sqrt{2}$. The choice $t=0$ implies $y_1=0,y_2=1,y_3=0$
which is inconsistent with the last but one constraint in
(\ref{schur3_opti}). The choice $t=1-\sqrt{2}$ yields
$y_1 = \pm (1+\sqrt{2})\alpha$, $y_2= 0$,
$y_3 = \pm \alpha$ with $\alpha=1/\sqrt{4-2\sqrt{2}}$
and the objective function $-2y_1^2+2y_1y_3 = -1+\sqrt{2}$.
The choice $t=1+\sqrt{2}$ yields $y_1 = \pm \alpha$, $y_2=0$,
$y_3 = \pm (-1-\sqrt{2})\alpha$ and the objective function
$-1-\sqrt{2}$, a negative minimum curvature. Therefore
region $S$ is indeed nonconvex.

From the remaining
constraints in (\ref{schur3_opti}), we conclude that
the minimal curvature points $x$ can be found along
the portion of parabola $\sqrt{2}x_1^2-x_2+1=0$  included
in the half-planes $(2+\sqrt{2})x_1+x_2+1 \geq 0$ and $-(2+\sqrt{2})x_1
+x_2+1 \geq 0$. Any plane tangent
to the hyperbolic paraboloid $p_3(x)=0$ at a point along
the parabola $\sqrt{2}x_1^2-x_2+1=0$ can be used to generate
a valid inner approximation of the stability region.
For example,
with the choice $x^1=(0,1,0)$, we generate the gradient
half-plane $p_4(x) = g_3(x^1)(x-x^1) = x_2 - 1 \leq 0$.

More generally, for discrete-time polynomials of degree $n\geq 3$,
stability region $S$ is the image of the box $B=[-1,1]^n$
(of so-called reflection coefficients) though a multiaffine
mapping, see e.g. \cite{n06} and references therein.
The boundary of $S$ consists of ruled surfaces, and
the convex hull of $S$ is generated by the images
of the vertices of $B$ through the multiaffine mapping.
It would be interesting to investigate whether this
particular geometry can be exploited to generate
systematically a convex inner approximation of maximum
volume of the stability region $S$.

\subsection{Fixed-order controller design}

Consider the open-loop discrete-time system $(-2z^2+1)/(z^3+z^2+a)$, parametrized
by $a \in \R$,
in negative closed-loop configuration with the controller $(x_1z+x_2)/(z+1)$.
The characteristic polynomial is equal to $q(z) = \sum_{k=0}^4 q_k z^k = z^4+2(1-x_1)z^3+(1-2x_2)z^2+(a+x_1)z+a+x_2$,
and it is Schur stable (all roots in the open unit disk) if and only if
$p_k<0$, $k=1,2,\ldots,4$
and $p_6=-p_2p_3p_4+p^2_2p_5+p_1p^2_4<0$ where
\[
\left[\begin{array}{c}
p_1 \\ p_2 \\ p_3 \\ p_4 \\ p_5
\end{array}\right] =
\left[\begin{array}{rrrrr}
-1 & 1 & -1 & 1 & -1 \\
4 & -2 & 0 & 2 & -4 \\
-6 & 0 & 2 & 0 & -6 \\
4 & 2 & 0 & -2 & -4 \\
-1 & -1 & -1 & -1 & -1
\end{array}\right]
\left[\begin{array}{c}
q_0 \\ q_1 \\ q_2 \\ q_3 \\ q_4 
\end{array}\right].
\]
The affine inequalities $p_k(x_1,x_2)<0$, $k=1,2,\ldots,5$ define a polytope
in the controller parameter plane $(x_1,x_2) \in \R^2$,
and the inequality $p_6(x_1,x_2)<0$ defines a cubic region.

\begin{figure}[htb]
\centering{\includegraphics[width=0.75\columnwidth]{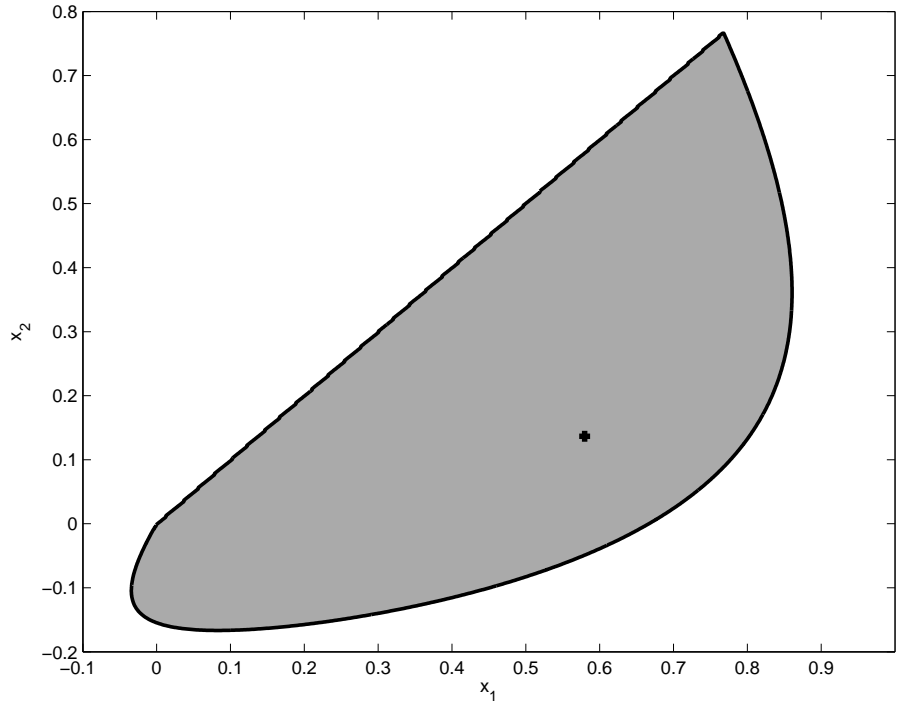}}
\caption{Convex stability region (dark gray),
with analytic center (cross) corresponding to a fixed-order controller.\label{degre4convex}}
\end{figure}

In the case $a=0$, with the following Gloptipoly 3 implementation of Steps 1-3 of Algorithm 1:
\begin{verbatim}
mpol x y 2
p1 = -x(1)+x(2);
p2 = -6*x(1)+4*x(2);
p3 = -10*x(2)-4;
p4 = -8+4*x(2)+6*x(1);
p5 = -4+x(1)+x(2);
p6 = 6*x(1)^2*x(2)+3*x(1)^2-10*x(1)*x(2)-2*x(1)-3*x(2)^3+6*x(2)^2+x(2);
g6 = diff(p6,x); % gradient
H6 = diff(g6,x); % Hessian
% LMI relaxation order
order = input('LMI relaxation order = ');
% build LMI relaxation
P = msdp(min(y'*H6*y), p6==0, p1<=0, p2<=0, p3<=0, p4<=0, p5<=0, ...
    g6*y==0, y'*y==1, order);
% solve LMI relaxation
[status,obj] = msol(P)
\end{verbatim}
we obtain a negative lower bound {\tt obj = -3.5583} at the 2nd LMI
relaxation, which is inconclusive. At the 3rd LMI relaxation,
we obtain a positive lower bound {\tt obj = 0.8973} which
certifies convexity of the stability region, see Figure \ref{degre4convex}.
Since the stability region $\bar{S} = \{x \in \R^2 \: :\: p_k(x) \leq 0, \:k=1,2,\ldots,6\}$
is convex, we can optimize over it with standard techniques of convex optimization. 
More specifically, a recent result in \cite{l11} indicates that any
limit point of any sequence of admissible stationary points
of the logarithmic barrier function $f(x) = -\sum_{k=1}^6 \log p_k(x)$
is a Karush-Kuhn-Tucker point satisfying first order optimality
condition. In particular, the gradient of $f(x)$
vanishes at the analytic center of the set. Using Maple
(or a numerical local optimization method) we can readily obtain
the analytic center $x^*_1 \approx 0.57975$, $x^*_2 \approx 0.13657$
(five-digit approximations of algebraic coefficients of degree 17)
corresponding to a controller well inside the stability region.
Such a controller can be considered as non-fragile, in the sense
that some uncertainty on its coefficients will not threaten 
closed-loop stability.

\begin{figure}[htb]
\centering{\includegraphics[width=0.75\columnwidth]{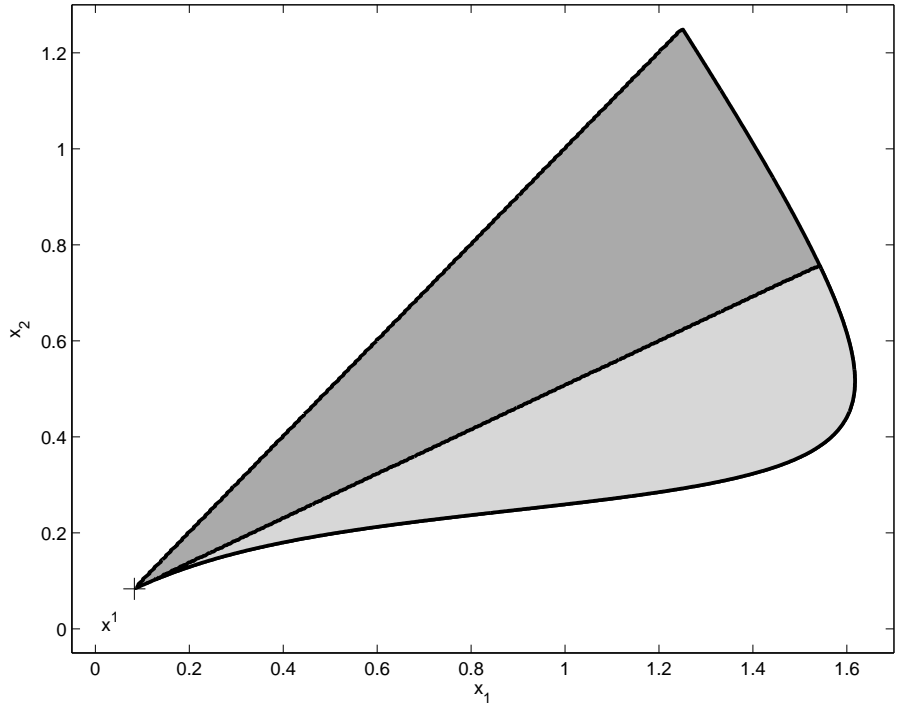}}
\caption{Convex inner approximation (dark gray) of nonconvex fourth-order
discrete-time stability region (light gray).\label{degre4nonconvex}}
\end{figure}

Now for the choice $a=-3/4$ we carry on again our study of convexity
of the stability region with the help of a similar GloptiPoly script.
At the 2nd LMI relaxation we obtain a negative lower bound {\tt obj = -385.14}
which is inconclusive. At the 3rd LMI relaxation, we obtain a negative
lower bound {\tt obj = -380.88} which is also inconclusive.
Eventually, at the 4th LMI relaxation, we obtain a negative
lower bound {\tt obj = -380.87} which is certified to
be the global minimum with {\tt status = 1}. The point $x^1$ at
which the minimum curvature is achieved is a vertex of the
stability region, and the tangent at this point
of the nonconvex part of the boundary is used to generate
a valid inner approximation $\bar{S}$, see Figure \ref{degre4nonconvex}.
Any point chosen in this triangular region corresponds to a stabilizing
controller.
We see that here the choice of the point of minimum curvature
is not optimal in terms of maximizing the surface of $\bar{S}$.
A point chosen elsewhere along the negatively curved part of
the boundary would be likely to generate a larger convex
inner approximation. 

\subsection{Optimal control with semialgebraic constraints}

In Model Predictive Control (MPC), an optimal control problem is solved recursively. This resolution is usually based
on direct methods that consist of deriving a nonlinear program from the optimal control problem by discretization of
the dynamics and the path constraints. Since the embedded software has strict specification on algorithm complexity and
realtime computation, convexity of the program is a key feature \cite{Ross}. Indeed, in this context,
our convex inner approximation of the admissible space become valuable to speed up the computation even
at the price of some conservatism.

In open-loop control design, convexity of the problem is a matter of concern especially when the optimal control problem
is part of an MPC procedure. In this case, the optimal control problem is solved mostly using direct
methods that transfom it into a parametric optimization problem. Convexity permits to limits
the complexity of the resolution and so reduces the computation time of an optimal solution. Unfortunately, in dynamic
inversion techniques based on differential flatness, the generally convex constraints on the states and inputs are replaced by
nonconvex admissible sets in the flat output space, see \cite{Ross} and reference therein for details. Thus, in such a
method, it is necessary to design inner convex approximation of the admissible subset to develop a tractable algorithm
\cite{Louembet09}.

\begin{figure}[htb]
\centering{\includegraphics[width=0.45\columnwidth]{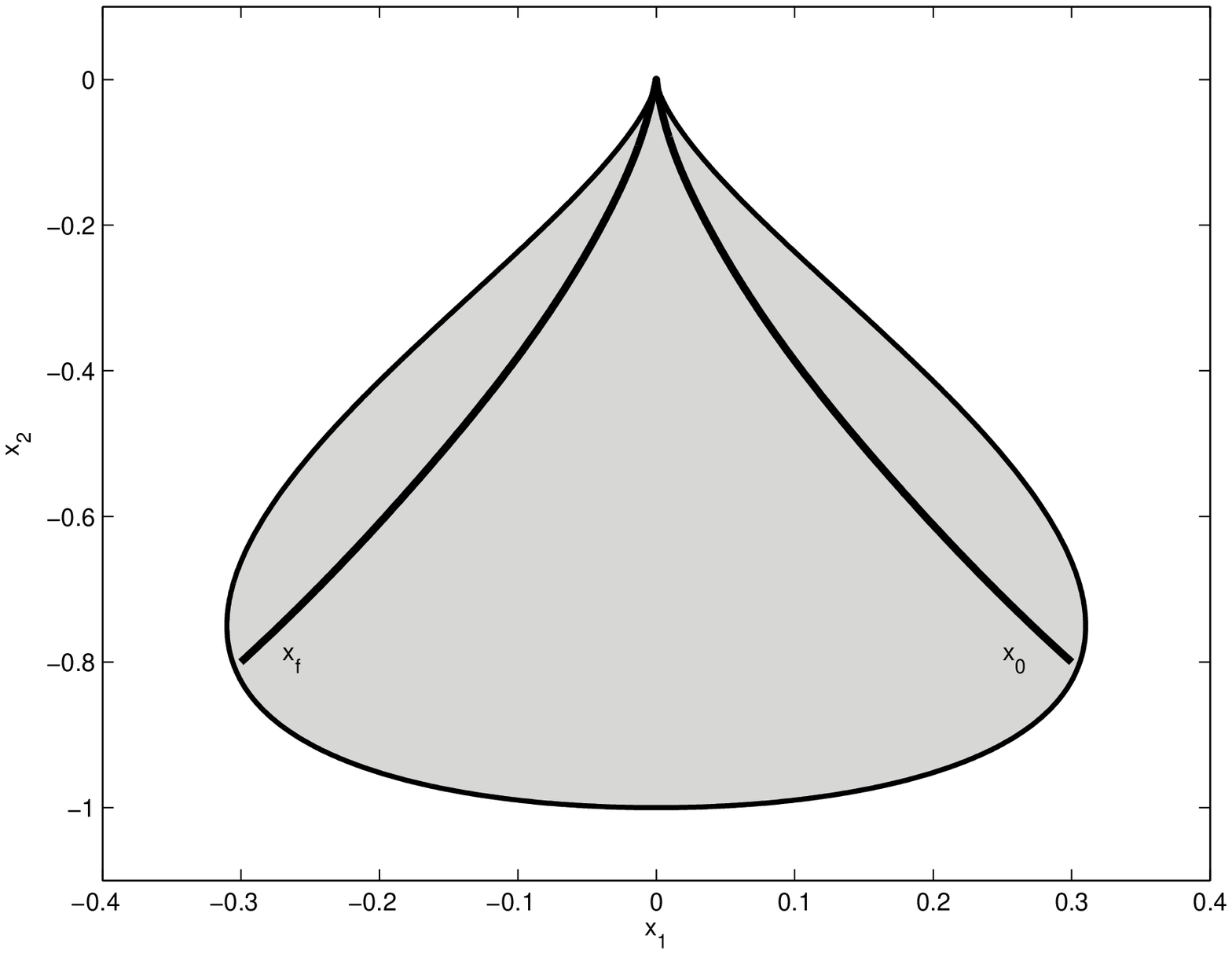}
\includegraphics[width=0.45\columnwidth]{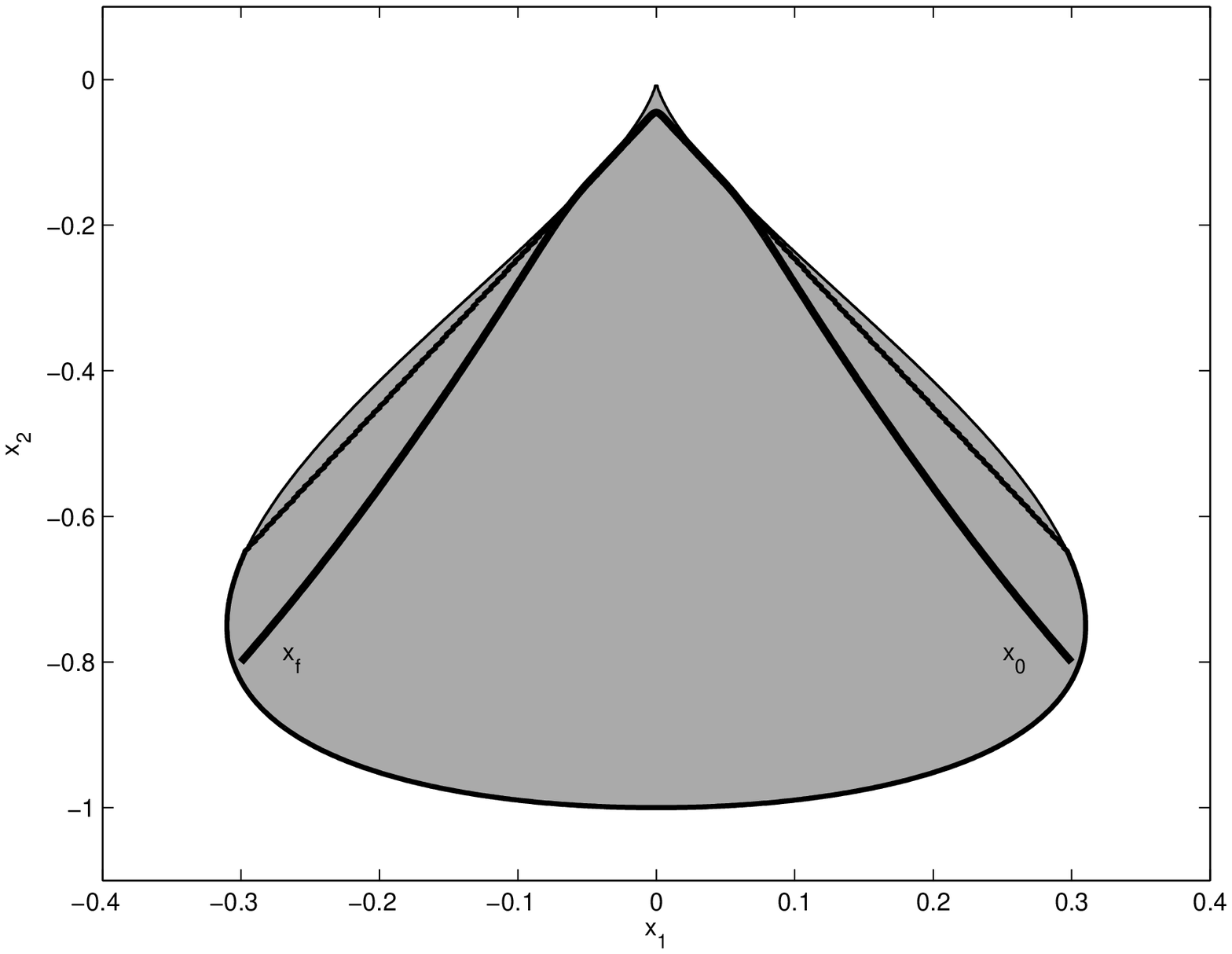}}
\caption{Optimal trajectories (bold) in nonconvex admissible set (left) and
in convex inner approximation (right).\label{waterdrop_optimalcontrol}}
\end{figure}

Consider the following optimal control problem
\[
\begin{array}{ll}
\min_{x,u} & \displaystyle \int_{t_0}^{t_f} u^2(t)dt \\
\text{s.t.} &
\dot{x} =\begin{bmatrix}0&1\\0&0\end{bmatrix}x+\begin{bmatrix}0\\1\end{bmatrix}u
\\
& x(t_0)=x_0, \quad x(t_f)=x_f\\
& p_1(x)\leq 0.
\end{array}
\]
The objective of this problem is to steer the linear system from an initial state to a final state in a fixed time inside
the admissible state subset $S$ defined e.g. by the waterdrop quartic defined in section \ref{secwaterdrop}:
\begin{equation}\label{contrainte_S}
S=\{x\in\mathbb R^2 \: :\: p_1(x)=x_1^4+x_2^4+x_1^2+x_2^3\leq 0\}.
\end{equation}
We describe thereafter a classical methodology for solving the previous optimal control problem using flatness-based
dynamic inversion, see \cite{Milam,Petit,Levine} for other examples.
As the dynamics are linear and fully actuated, dynamic inversion can be used to develop an efficient algorithm for
the considered problem \cite{Levine}. Thus, the system trajectory can be parametrized by a user-specified sufficiently smooth
function $x_1(t)=f(t)$
so that $x_2(t)=\dot f(t)$ and $u(t)=\ddot f(t)$. The function $f(t)$ is classically described by a chosen basis $b(t)$ and
the associated vector of weighting coefficient $\alpha$ such that
\[f(t)=\sum_k \alpha_k b_k(t).\] 
In order to derive a finite dimensional program, the admissible set constraint is discretized and
enforced at a finite number of time instants $\{t_i\}_{i=1,\dots,N}$ such that $t_0\leq t_1<t_2<\dots<t_N\leq t_f$.
Since $\mathcal S$ is nonconvex, we obtain a finite-dimensional nonlinear nonconvex program:
\[
\begin{array}{ll}
\min_{\alpha} & \sum_i u^2(\alpha,t_i) \\
\text{s.t.} & x(\alpha,t_0)=x_0, \quad x(\alpha,t_f)=x_f\\
& p_1(x(\alpha,t_i))\leq 0, \quad i=1,\dots,N.
\end{array}
\]
The inner approximation $\bar S$ calculated previously in section \ref{secwaterdrop} is given by
$\bar{S}=\{x\in\mathbb R^2 \: :\: p_1(x)\leq 0,~p_2(x)=g_1(x^1)(x-x^1)\leq 0,~p_3(x)=g_1(x^2)(x-x^2)\leq 0 \}$
where $g_1(x^1)$ and $g_1(x^2)$ is the gradient of $p_1(x)$ evaluated at $x=x^1$ and $x=x^2$, respectively.
The use of the inner approximation $\bar S$ as admissible subset leads to the following convex program:
\[
\begin{array}{ll}
\min_{\alpha} & \sum_i u^2(\alpha,t_i)\\
\text{s.t.} & x(\alpha,t_0)=x_0, \quad x(\alpha,t_f)=x_f\\
& p_1(x(\alpha,t_i))\leq 0, \quad p_2(x(\alpha,t_i))\leq 0, \quad p_3(x(\alpha,t_i))\leq 0, \quad i=1,\dots,N.
\end{array}
\]
In the following, we set $t_0=0$, $x_0=[0.3000,\:-0.8000]$ and $t_f=2.5$, $x_f=[-0.3000,\:-0.8000]$.
The time function $f(t)$ is a 5-segment-piecewise polynomial of the
4th order (degree 3) defined on a B-spline basis. We run both programs for different values of $N$. In table \ref{table}
we compare the computation times and optimal costs. See Figure \ref{waterdrop_optimalcontrol} for
the state trajectories.
\begin{table}[h!]
\begin{center}
\begin{tabular}{c|c|cccccccc}
& $N$ &10& 20& 50& 100 & 200& 500& 1000\\\hline
CPU time [s] & Convex& 0.029& 0.045& 0.056& 0.058&0.102 & 0.225 & 0.498\\\hline
        & Nonconvex& 0.109&0.195 &0.199 &0.409 &0.513 &0.836 & 1.46\\\hline
Optimal cost & Convex& 1.65& 1.67&1.68 &1.69 &1.68 & 1.68& 1.68 \\\hline
        & Nonconvex&1.52 &1.52 &1.52 &1.52 &1.52 &1.52 &1.52\\\hline
\end{tabular}
\caption{Computation times and optimal costs of
the nonconvex and convexified optimal control problems,
as functions of the number $N$ of discretization points.\label{table}}
\end{center}
\end{table}
For this example, we observe the positive effect that convexity has on the
reduction of the computational burden, balanced by the relatively small
loss of performance.

\section{Conclusion}

We have presented a general-purpose computational algorithm to
generate a convex inner approximation of a given basic semialgebraic
set. The inner approximation is not guaranteed to be of maximum
volume, but the algorithm has the favorable features of
leaving invariant a convex set, and preserving convex boundaries
while removing nonconvex regions by enforcing linear constraints
at points of minimum curvature. Even though our initial motivation
was to construct convex inner approximations of stability regions
for fixed-order controller design, our algorithm can be used
on its own for checking convexity of semialgebraic sets.

Each step of the algorithm consists
in solving a potentially nonconvex polynomial optimization problem
with the help of a hierarchy of convex LMI relaxations. For this
we use Gloptipoly 3, unfortunately with no guarantee of a priori
computational burden, even though in practice it is observed
that global optimality is ensured at a moderate cost, as soon
as the dimension of the ambient space is small. Numerical
experiments indicate that the approach may be practical for
ambient dimensions up to 4 or 5. For larger problems, we can
rely on more sophisticated nonlinear or global optimization codes
\cite{n04}, even though this possibility has not been investigated
in this paper. Indeed, our main driving force is to contribute with
a readily available Matlab implementation.

Our algorithm returns a sequence of polynomials such that the
intersection  of their sublevel sets is geometrically convex.
However, the individual polynomials (of degree two or more) are not
necessarily convex functions. One may therefore question the relevance
of applying a relatively complex algorithm to obtain a
convex inner approximation in the form of a list of defining
polynomials which are not necessary individually convex.
A recent result of \cite{l11} indicates however that
any local optimization method based on standard first-order
optimality conditions for logarithmic barrier functions
will generate a sequence of iterates converging to
the global minimum of a convex function over
convex sets. In other words,
geometric convexity seems to be more important
that convexity of the individual defining polynomials. 

Indeed, if convexity of the inner approximation is guaranteed
in the presented work, convexity of the defining polynomials
would allow the use of constant multipliers to certificate
optimality in a nonlinear optimization framework.
Instead, with no guarantee of convexity of the defining polynomials,
the geometric proprety of convexity of the sets is more delicate
to exploit efficiently by optimization algorithms.

Finally, let us emphasized that it is conjectured that all convex semialgebraic sets
are semidefinite representable in \cite{hn07}, see also \cite{l06}.
It may then become possible to fully exploit the geometric
convexity of our inner convex through an explicit representation
as a projection of an affine section of the semidefinite cone.
For example, in our target application domain, this would allow
to use semidefinite programming to find a suboptimal stabilizing
fixed-order controller.

\section*{Acknowledgements}

The first author is grateful to J. W. Helton for many discussions
and ideas leading to Theorem \ref{convexity_thm}.

\end{document}